\documentclass[11pt,reqno]{amsart}

\usepackage[ams]{optional}

\usepackage{amssymb,amsopn,amsmath,tikz}
\usepackage[colorlinks]{hyperref}
\usepackage[all]{xy}

\opt{lms}{\title[Moduli spaces of torsion sheaves on K3 surfaces]{Moduli spaces of torsion sheaves on K3~surfaces and derived equivalences}}
\opt{ams}{\title[Moduli spaces of torsion sheaves on K3 surfaces\dots]{Moduli spaces of torsion sheaves on K3~surfaces and derived equivalences}}

\opt{lms}{
\author{Nicolas Addington, Will Donovan and Ciaran Meachan}
}

\opt{ams}{
\author[N.~Addington]{Nicolas Addington}
\address{Nicolas Addington \\
Department of Mathematics \\
University of Oregon \\
Eugene, OR 97403-1222 \\
United States}
\email{adding@uoregon.edu}

\author[W.~Donovan]{Will Donovan}
\address{Will Donovan \\
Kavli Institute for the Physics and Mathematics of the Universe (WPI) \\
The University of Tokyo Institutes for Advanced Study \\
5-1-5 Kashiwanoha \\
Kashiwa, Chiba, 277-8583 \\
Japan}
\email{will.donovan@ipmu.jp}

\author[C.~Meachan]{Ciaran Meachan}
\address{Ciaran Meachan \\ School of Mathematics \\
The University of Edinburgh \\
James Clerk Maxwell Building \\
Peter Guthrie Tait Road \\
Edinburgh EH9 3FD \\
United Kingdom}
\email{ciaran.meachan@ed.ac.uk}
}

\newcommand \C {\mathcal C}
\newcommand \CC {\mathbb C}
\newcommand \E {\mathcal E}
\newcommand \F {\mathcal F}
\renewcommand \H {\mathcal H}
\newcommand \I {\mathcal I}
\renewcommand \L {\mathcal L}
\newcommand \M {\mathcal M}
\newcommand \N {\mathcal N}
\renewcommand \O {\mathcal O}
\renewcommand \P {\mathbb P}
\newcommand \U {\mathcal U}
\newcommand \Z {\mathbb Z}
\renewcommand \a {\mathbf a}
\newcommand \s {\mathbf s}
\renewcommand \v {\mathbf v}

\DeclareMathOperator \Br {Br}
\DeclareMathOperator \Hilb {Hilb}
\newcommand \Picbar {\overline{\mathrm{Pic}}{}} % extra {} to keep any following superscript from getting too high
\newcommand \FM {F\!M}
\newcommand \KN {K\!N}
\DeclareMathOperator \id {id}
\DeclareMathOperator \Ext {Ext}
\DeclareMathOperator \cone {cone}
\DeclareMathOperator \RHom {RHom}
\DeclareMathOperator \Pic {Pic}
\DeclareMathOperator \RGamma {R\Gamma}
\DeclareMathOperator \rank {rank}
\DeclareMathOperator \Hom {Hom}
\DeclareMathOperator \NS {NS}
\DeclareMathOperator \Mov {Mov}
\DeclareMathOperator \Nef {Nef}

\newcommand \Caldararu {C\u{a}ld\u{a}\-raru}

\opt{lms}{}
\opt{ams}{\newtheorem*{conj*}{Conjecture}}
\newtheorem{thm}{Theorem}[section]
\newtheorem{prop}[thm]{Proposition}

\newtheorem{lem}[thm]{Lemma}
\opt{ams}{\theoremstyle{definition}}
\newtheorem{rmk}[thm]{Remark}
\newtheorem{defn}[thm]{Definition}

\opt{ams}{
\numberwithin{equation}{section}

}

\begin{document}

% !TeX root = main.tex

% 14F05  	Sheaves, derived categories of sheaves and related constructions
% 14J60   	Vector bundles on surfaces and higher-dimensional varieties, and their moduli
% 14J28  	$K3$ surfaces and Enriques surfaces
% 14K30   	Picard schemes, higher Jacobians
% 14E30   	Minimal model program (Mori theory, extremal rays)
% 53C26   	Hyper-Kähler and quaternionic Kähler geometry, "special'' geometry

\opt{lms}{
\classno{Primary 14F05; Secondary 14J60, 14J28, 14K30, 14E30, 53C26 }
\extraline{The first author was partially supported by NSF grant no.\ DMS--0905923.  The second author was supported by World Premier International Research Center Initiative (WPI Initiative), MEXT, Japan; and by EPSRC grant~EP/G007632/1.  The third author was supported by the EPSRC Doctoral Prize Research Fellowship Grant EP/K503034/1.}
}

\opt{lms}{\maketitle}

\begin{abstract}
We show that for many moduli spaces $\M$ of torsion sheaves on K3 surfaces $S$, the functor $D^b(S) \to D^b(\M)$ induced by the universal sheaf is a $\P$-functor, hence can be used to construct an autoequivalence of $D^b(\M)$, and that this autoequivalence can be factored into geometrically meaningful equivalences associated to abelian fibrations and Mukai flops.  Along the way we produce a derived equivalence between two compact hyperk\"ahler $2g$-folds that are not birational, for every $g \ge 2$.  We also speculate about an approach to showing that birational moduli spaces of sheaves on K3 surfaces are derived-equivalent.
\end{abstract}

\opt{ams}{\maketitle}

\section*{Introduction}

The group of autoequivalences of the derived category of coherent sheaves on a variety is an interesting and subtle geometric object.  Any autoequivalences beyond the ``standard'' ones -- automorphisms of the variety itself, tensoring by line bundles, and homological shift -- should be seen as ``hidden symmetries'' of the variety \cite{bridgeland_icm}.  Many non-standard autoequivalences come from birational geometry: twists of one kind or another around subvarieties that can be contracted or flopped.  In \cite{nick} the first author introduced a rather different autoequivalence for the Hilbert scheme of $n$ points on a K3 surface, built from the universal ideal sheaf, and conjectured that the same construction would work for any moduli space of sheaves on a K3 surface.  Our main goal in this paper is to prove this conjecture for certain moduli spaces of torsion sheaves:
{ \renewcommand{\thethm}{A}
\begin{thm} \label{main_thm}
Let $S$ be a complex projective K3 surface of Picard rank 1 and degree $2g-2$, and let $\M$ be the moduli space of stable sheaves with Mukai vector $(0,1,d+1-g)$.  Let $\alpha \in \Br(\M)$ be the Brauer class obstructing the existence of a universal sheaf, and let $\F$ be a $(1 \boxtimes \alpha)$-twisted pseudo-universal sheaf on $S \times \M$.  Then the functor $F\colon D^b(S) \to D^b(\M,\alpha)$ induced by $\F$ is a $\P^{g-1}$-functor, hence can be used to construct an autoequivalence of $D^b(\M,\alpha)$.
\end{thm} }
\noindent The general member of $\M$ is a degree-$d$ line bundle on a smooth genus-$g$ curve in $S$, so $\M$ is a $2g$-fold fibered over $\P^g$ in Jacobians.  If $2g-2$ is relatively prime to $d+1-g$ then the Brauer class $\alpha$ vanishes and $\F$ is an honest sheaf.  We briefly review the definition of $\P$-functors and the associated autoequivalences, called $\P$-twists, in \S\ref{P-review}.

Whereas the geometric meaning of the earlier autoequivalence of \linebreak $D^b(\Hilb^n(S))$ was somewhat obscure, our new autoequivalence of $D^b(\M)$ turns out to factor as a product of Fourier--Mukai--Arinkin equivalences associated to abelian fibrations and Kawamata--Namikawa equivalences associated to Mukai flops.  Moreover, when $n=g$ the old autoequivalence is conjugate to the new one by a Kawamata--Namikawa equivalence.  We review these equivalences and their relation to $\P$-twists in \S\ref{KN-review}.

\subsection*{Outline of the argument}  Assume the set-up of Theorem \ref{main_thm}, let $\C$ be the universal curve over the linear system $|\O_S(1)| = \P^g$, and identify $\M$ with the compactified relative Picard variety $\Picbar^d := \Picbar^d(\C/\P^g)$.  In \S\ref{fma} we show that $F\colon D^b(S) \to D^b(\Picbar^d, \alpha)$ factors as
\begin{equation} \label{factor_F}
F = \FM \circ AJ_* \circ \varpi^*(- \otimes \O_S(l)),
\end{equation}
where
\setlength \leftmargini {1.5em}
\begin{itemize}
\item $l$ is an integer,
\item $\varpi$ is the natural map $\C \to S$, which is a $\P^{g-1}$-bundle,
\item $AJ\colon \C \to \Picbar^{-1}$ is the Abel--Jacobi embedding, and
\item $\FM \colon D^b(\Picbar^{-1}, \beta) \to D^b(\Picbar^d, \alpha)$ is a family version of Mukai's derived equivalence between an abelian variety and its dual \cite{fourier_mukai}.  The extension to the singular fibers is due to Arinkin \cite{arinkin}.  Here $\beta$ is a Brauer class on $\Picbar^{-1}$ with $AJ^*(\beta) = 1$.
\end{itemize}
\setlength \leftmargini {2em}
Now since $\FM$ and $- \otimes \O_S(l)$ are equivalences, it remains to show that $AJ_* \circ \varpi^*$ is a $\P^{g-1}$-functor, and to understand the associated $\P$-twist.  Since $AJ(\C) \subset \Picbar^{-1}$ is the center of a Mukai flop $\Picbar^{-1} \dashrightarrow X'$, we can use our results from \cite{flop_flop}; in Proposition \ref{flop_flop_with_brauer_classes} we check that the Brauer class $\beta$ does not cause trouble.  The upshot is that the $\P$-twist $P_F \in \operatorname{Aut}(D^b(\Picbar^d,\alpha))$ factors as
\begin{equation} \label{factor_PPic}
P_F = \FM \circ \KN^{-1}_{g-1} \circ \KN_g \circ \FM^{-1},
\end{equation}
where
\[ \KN_k \colon D^b(\Picbar^{-1},\beta) \to D^b(X',\beta') \qquad \qquad k\in \Z \]
are Kawamata--Namikawa equivalences associated to the Mukai flop.
\bigskip

Along the way we obtain the following result of independent interest:
{ \renewcommand{\thethm}{B}
\begin{thm} \label{per_not_inv}
For every $g \ge 2$ there are compact hyperk\"ahler $2g$-folds $X$ and $Y$ such that $D^b(X) \cong D^b(Y)$ but $H^2(X,\Z)$ is not Hodge-isometric to $H^2(Y,\Z)$.  In particular $X$ is not birational to $Y$.
\end{thm} }
\noindent Precisely, we take $X = \Picbar^0$ and $Y = \Picbar^{g-1}$.

\subsection*{Intertwinement with the Hilbert scheme}
There is a beautiful bi\-rational map $\Hilb^g(S) \dashrightarrow \Picbar^{-g}$, defined as follows: a generic set of $g$ points $\zeta \subset S$ is contained in a unique curve $C \in |\O_S(1)|$, and we send $[\zeta] \in \Hilb^g$ to $[\O_C(-\zeta)] \in \Picbar^{-g}$.  For $g=2$, this is the original example of a Mukai flop \cite[Example 0.6]{mukai_inventiones}: we have a double cover $f\colon S \to \P^2$ branched over a smooth sextic curve; if the two points of $\zeta$ map to distinct points in $\P^2$ then $C$ is the preimage in $S$ of the line that they span in $\P^2$; the indeterminacy locus of the flop is the Lagrangian $\P^2$ in $\Hilb^2(S)$ consisting of length-2 subschemes of the form $f^{-1}(\text{point})$.  The details are equally pretty for $g=3$, which also appeared in Mukai's original paper \cite[Example 0.8]{mukai_inventiones}.

On the level of sheaves, this birational map $\Hilb^g \dashrightarrow \Picbar^{-g}$ is implemented by the spherical twist around $\O_S(-1)$: if $\zeta \subset S$ is contained in a unique curve $C$ then the twist sends $\I_{\zeta/S}$ to $\O_C(-\zeta)$.  (On the indeterminacy locus, where $\zeta$ is contained in a pencil of curves, the twist sends $\I_{\zeta/S}$ to a two-term complex, not a sheaf.)  Thus one is led to look for some compatibility between the functors
\begin{align*}
F \colon D^b(S) &\to D^b(\Picbar^{-g}) &
F' \colon D^b(S) &\to D^b(\Hilb^g)
\end{align*}
induced by the universal sheaves.  In \S\ref{compatibility} we show that if $g \le 5$, so the birational map is a Mukai flop, then
\begin{equation} \label{compat_eq}
F' \circ T_{\O_S(1)} = \KN_2 \circ F,
\end{equation}
where $\KN_2\colon D^b(\Picbar^{-g}) \to D^b(\Hilb^g)$ is again a Kawamata--Namikawa equivalence associated to the Mukai flop.  As a consequence, the associated $\P$-twists are conjugate:
\[ P_{F'} = \KN_2 \circ P_F \circ \KN_2^{-1}. \]
Together with the factorization \eqref{factor_PPic} this gives
\[ P_{F'} = \KN_2 \circ \FM \circ \KN^{-1}_{1-g} \circ \KN_g \circ \FM^{-1} \circ \KN_2^{-1}, \]
so the autoequivalence of $D^b(\Hilb^g)$ factors into a product of geometrically meaningful equivalences as promised.

For $g > 5$, the birational map $\Hilb^g \dashrightarrow \Picbar^{-g}$ is not a Mukai flop but a \emph{stratified Mukai flop} \cite{brill_noether}.  The relevant equivalences are due to Cautis, Kamnitzer, and Licata \cite{ckl,cautis_stratified}, and we expect that they can be adapted to our situation to give a formula like \eqref{compat_eq}, but we do not pursue this as the bookkeeping begins to overwhelm the geometry.

\subsection*{Approach to flops between moduli spaces in general}  We conclude this introduction with some speculation about how the compatibility \eqref{compat_eq} might fit into a broader framework.  Let $S$ be any K3 surface and $\M$ any fine moduli space of stable sheaves, or more generally of $\sigma$-stable objects for some Bridgeland stability condition $\sigma$.  A well-known conjecture of Bondal and Orlov \cite[Conj.~4.4]{bo_icm} and Kawamata \cite[Conj.~1.2]{kawamata} would imply that any other smooth $K$-trivial birational model $\M'$ of $\M$ has $D^b(\M') \cong D^b(\M)$.  Recently Bayer and Macr\`i \cite[Thm.~1.2(c)]{bm13} have shown that for any such $\M'$ there is another stability condition $\sigma'$ such that $\M'$ is the moduli space of $\sigma'$-stable objects with the same Mukai vector.  We propose the following seemingly stronger conjecture:
\begin{conj*}
With $\M$ and $\M'$ as in the previous paragraph, there is an equivalence $\Phi\colon D^b(\M) \to D^b(\M')$ such that for suitably normalized universal objects $\F \in D^b(S \times \M)$ and $\F' \in D^b(S \times \M')$, the induced functors $F\colon D^b(S) \to D^b(\M)$ and $F'\colon D^b(S) \to D^b(\M')$ satisfy $F' = \Phi \circ F$.  In other words, the equivalence ${\id} \boxtimes \Phi\colon D^b(S \times \M) \to D^b(S \times \M')$ takes $\F$ to $\F'$.
\end{conj*}

Our compatibility \eqref{compat_eq} verifies this conjecture for $\M = \Hilb^g$ and $\M' = \Picbar^{-g}$; \emph{a priori} these are moduli spaces for the same stability condition and different Mukai vectors, but pulling back the stability condition and the second Mukai vector via $T_{\O_S(-1)}$ puts us in the situation of the conjecture. For good measure we show in \S\ref{movable} that $\Hilb^g$ has no other $K$-trivial birational models apart from $\Picbar^{-g}$.

The hard way to go about proving the conjecture in general would be to produce an equivalence $\Phi$, say by some combinatorial recipe from a stratification of the indeterminacy locus of the birational map $\M \dashrightarrow \M'$, and then go on to check that it takes $\F$ to $\F'$.  But Markman has suggested that one might instead try to construct $\Phi$ \emph{from} $\F$ and $\F'$.  He showed in \cite[Thm~1.2(1)]{on_the_monodromy} that the middle Chern class of the kernel $F' \circ R$, where $R\colon D^b(\M) \to D^b(S)$ is the right adjoint of $F$, induces an isomorphism of cohomology rings $H^*(\M,\Z) \cong H^*(\M',\Z)$.  One has $\H^0(F \circ R) = \O_\Delta$; perhaps in general $\H^0(F' \circ R)$ is the kernel of an equivalence?

\opt{lms}{\begin{acknowledgements}}\opt{ams}{\subsection*{Acknowledgements}}
We thank Eyal Markman for several good ideas, % in Kyoto, suggested that \KN gives compatibility between Hilb^g and Pic^g; in Lille, suggested Bayer-Macri speculation at end of intro.
Justin Sawon % pointed us to Arinkin paper
and Dima Arinkin % helped us sort out Sawon issue
for expert advice, and Arend Bayer for helpful discussions. % Ciaran talked to him a lot
Much of this work was done while visiting the Hausdorff Research Institute for Mathematics in Bonn, Germany, and we thank them for their hospitality.  
\opt{ams}{N.A.\ was partially supported by NSF grant no.\ DMS--0905923.  W.D.\ was supported by World Premier International Research Center Initiative (WPI Initiative), MEXT, Japan; and by EPSRC grant~EP/G007632/1.  C.M.\ was supported by the EPSRC Doctoral Prize Research Fellowship Grant EP/K503034/1.}
\opt{lms}{\end{acknowledgements}}
% !TeX root = main.tex

\section{Review of \texorpdfstring{$\P$}{P}-functors and Mukai flops}

\subsection{\texorpdfstring{$\P$}{P}-functors} \label{P-review}
A \emph{$\P^n$-object} on a hyperk\"ahler $2n$-fold $X$ is an object $\E \in D^b(X)$ such that $\Ext^*_X(\E,\E) \cong H^*(\P^n, \CC)$ as rings.  The main examples are skyscraper sheaves of Lagrangian $\P^n$s in $X$, if there are any, and line bundles on $X$.  From such an object, Huybrechts and Thomas \cite{ht} constructed an autoequivalence
\[ P_\E\colon D^b(X) \to D^b(X), \]
as a certain double cone
\[ P_\E(\F) = \cone\!\big(\cone\!\big(\E \otimes \RHom(\E,\F)[-2] \to \E \otimes \RHom(\E,\F)\big) \xrightarrow{\text{eval}} \F\;\big) \]
This equivalence is called the \emph{$\P$-twist around $\E$}.

In \cite{nick}, the first author observed that if $S$ is a K3 surface and $\Hilb^n = \Hilb^n(S)$ is the Hilbert scheme of $n$ points on $S$, then the universal ideal sheaf $\I$ on $S \times \Hilb^n$ is in some sense a relative $\P^{n-1}$-object over $S$.  One might expect this to mean that $\E xt^*_S(\I,\I) \cong \O_S \otimes H^*(\P^n, \CC)$, but this is not true and, luckily, it is not what is needed in order to generalize Huybrechts and Thomas's twists.  Instead, consider the functor $F\colon D^b(S) \to D^b(\Hilb^n)$ induced by $\I$.  Then for $\F_1, \F_2 \in D^b(S)$ we have
\begin{equation} \label{ext-version}
\Ext^*_{\Hilb^n}(F(\F_1), F(\F_2)) \cong \Ext^*_S(\F_1, \F_2) \otimes H^*(\P^{n-1}, \CC)
\end{equation}
``as rings'': that is, if we have three objects $\F_1, \F_2, \F_3 \in D^b(S)$ then composition on the left-hand side of \eqref{ext-version} agrees with composition in the first factor on the right-hand side and the ring structure in the second factor.  A more compact way of saying this is to let $R \colon D^b(\Hilb^n) \to D^b(S)$ be the right adjoint of $F$ and to require that
\[ RF \cong \id_S \otimes H^*(\P^{n-1}, \CC) \]
and that the monad structure $RFRF \xrightarrow{R\epsilon F} RF$ agrees with the ring structure on $H^*(\P^{n-1}, \CC)$, at least on the level of cohomology sheaves.  For a more careful and more general definition of $\P$-functors, and a fuller discussion, see \cite[\S3]{nick} or \cite[\S1]{flop_flop}.

From a $\P$-functor $F\colon D^b(Z) \to D^b(X)$ the first author then constructed an autoequivalence
\[ P_F\colon D^b(X) \to D^b(X) \]
as a certain double cone
\[ P_F = \cone\!\big(\cone(FR[-2] \to FR) \xrightarrow\epsilon \id_X\big). \]
which reduces to Huybrechts and Thomas's construction when $Z$ is a point.  We do not need the precise definition here, but only some key facts:
\begin{enumerate}
\item If all the objects in the image of $F$ are supported on a subvariety $Y \subset X$, then $P_F$ acts as the identity on objects supported on $X \setminus Y$ and as a shift by $-2n$ on $\operatorname{im}(F)$.  Hence it is really a non-standard autoequivalence, not a composition of automorphisms of $X$ and line bundles and shifts.
\item If $\Phi\colon D^b(X) \to D^b(X')$ is an equivalence then $P_{\Phi F} \cong \Phi P_F \Phi^{-1}$.
\item \label{third_fact} If $\Psi\colon D^b(Z') \to D^b(Z)$ is an equivalence then $P_{F \Psi} \cong P_F$.
\end{enumerate}
\bigskip

Theorem \ref{main_thm} and the result on Hilbert schemes just described are instances of the following conjecture.\footnote{Markman and Mehrotra \cite[Cor.~6.18]{mm} have just proved another case of this conjecture, where the Brauer class $\alpha$ has order $2n-2$, the maximum possible.  Their proof is conditional on a certain conjecture about hyperholomorphic sheaves.}
\begin{conj*}[\cite{nick}]
Let $\M$ be a smooth, compact, $2n$-dimensional moduli space of sheaves on a K3 surface $S$, let $\alpha \in \Br(\M)$ be the Brauer class obstructing the existence of a universal sheaf, and let $\U$ be a $(1 \boxtimes \alpha)$-twisted pseudo-universal sheaf on $S \times \M$.  Then the functor $D^b(S) \to D^b(\M,\alpha)$ induced by $\U$ is a $\P^{n-1}$-functor.
\end{conj*}

The difficulty in proving this conjecture is that one needs to view points $x \in S$ as parametrizing (twisted) sheaves $\U|_{x \times \M}$ on $\M$, and understand the Ext groups between them; whereas one is accustomed to viewing points of $\M$ as parametrizing sheaves on $S$.  With Hilbert schemes and with moduli spaces of torsion sheaves (as in Theorem \ref{main_thm}) we know enough about the geometry of $\M$ and $\U$ to come to grips with these ``wrong-way slices'' of $\U$, but in general we know very little about them.  One approach to proving the conjecture in general would be to start with an elliptic K3 surface and induct from the Hilbert scheme case to higher-rank sheaves following O'Grady's careful analysis \cite{og}, which Marian and Oprea have interpreted in Fourier--Mukai terms in \cite[\S2.2]{mo}.  From there one could deform to a generic K3 surface easily enough; but getting to \emph{all} K3 surfaces would require major progress in understanding stability conditions in higher dimensions, or else some other new idea.
\bigskip

Other examples of $\P$-functors have been given by the third author \cite{ciaran_kummer} for generalized Kummer varieties, and by Krug \cite{krug_diagonal,krug_nakajima} for Hilbert schemes of other surfaces, using objects supported on certain correspondences rather than the universal ideal sheaf.  

In \cite{flop_flop} we studied $\P$-twists associated to Lagrangian $\P^n$s and coisotropic $\P^n$-bundles, in connection with Mukai flops, which we now review.

\subsection{Mukai flops and Kawamata--Namikawa equivalences} \label{KN-review}

Suppose that we have
\[ \xymatrix{
P \ar[d]_\varpi \ar[r]^j & X \\
{\phantom,}B,
} \]
where $X$ is a projective hyperk\"ahler variety, $P$ is a $\P^n$-bundle over a smooth, projective, connected base $B$, and $j$ is a closed embedding of codimension $n$.  Then the normal bundle $\N_{P/X}$ is isomorphic to the relative cotangent bundle $\Omega^1_{P/B}$ by \cite[Prop.~3.1(2)]{mukai_inventiones}, and we can consider the Mukai flop of $X$ along $P$:
\[ \xymatrix{
& \tilde X \ar[ld]_q \ar[rd]^p \\
X & & X'.
} \]
Here $\tilde X$ is the blow-up of $X$ along $P$, or of $X'$ along the dual $\P^n$-bundle $P^*$.  We assume that $X'$ is also projective.

Let $E \subset \tilde X$ be the exceptional divisor, which is identified with the universal hyperplane in $P \times_B P^*$, and let
\[ \hat X = \tilde X \cup_E (P \times_B P^*). \]
The line bundles $\O_{\tilde X}(E)$ and $\O_{P \times_B P^*}(-E)$ become isomorphic when restricted to $E$, and there is a unique way to glue them to get a line bundle on $\hat X$, which we call $\L$.

\begin{defn} \label{KN_k}
Let $X$, $X'$, and $\L \in \Pic(\hat X)$ be as in the previous paragraph.  For $k \in \Z$, we define
\[ D^b(X) \xrightarrow{\quad \KN_k \quad} D^b(X') \]
to be the functor induced by $\L^{\otimes k}$.
\end{defn}
\noindent Kawamata \cite[\S5]{kawamata} and Namikawa \cite{namikawa} showed that $\KN_0$ is an equivalence; for a textbook account see \cite[\S11.4]{huybrechts_fm}.  The extension to arbitrary $k$ is straightforward.

Assume that $P$ is the projectivization of a vector bundle, so we can speak about $\O_{P/B}(k)$ for all $k \in \Z$.  In \cite{flop_flop} we proved the following, which generalizes an example due to Cautis \cite[Prop.~6.8]{cautis_about}:

\begin{thm}[{\cite[Thm.~B$'$]{flop_flop}}] \label{flop_flop_thm}
If $H\!H^\text{odd}(B) = 0$ then the functor
\[ F_k = j_*(\O_{P/B}(k) \otimes \varpi^*(-)) \colon D^b(B) \to D^b(X) \]
is a $\P^n$-functor, and the associated $\P^n$-twist satisfies
\begin{equation} \label{flop_flop_fmla}
P_k \cong \KN_{n+k}^{-1} \circ \KN_{n+k+1}.
\end{equation}
\end{thm}

For our application we need a small extension of this theorem, to the case where $X$ carries a Brauer class $\beta$ with $j^* \beta = 1$.  The Brauer group is a birational invariant of smooth projective varieties, so $\beta$ determines a class $\beta' \in \Br(X')$, and we want to say there are equivalences $\KN_k\colon D^b(X,\beta) \to D^b(X',\beta')$ satisfying a version of Theorem \ref{flop_flop_thm}.  There is one subtlety, in that in order to define $D^b(X,\beta)$ and $D^b(X',\beta')$ we must choose cocycles representing $\beta$ and $\beta'$, and some care is required to make the choices compatibly.  We do this under a mild hypothesis:

\begin{prop} \label{flop_flop_with_brauer_classes}
Let $\beta \in \Br(X)$ and $\beta' \in \Br(X')$ be as in the previous paragraph, and suppose that there is a line bundle on $X$ whose restriction to $P$ is $\O_{P/B}(1)$.  Then for a suitable choice of cocycles representing $\beta$ and $\beta'$, the equivalences
\[ D^b(X,\beta) \xrightarrow{\quad \KN_k \quad} D^b(X',\beta') \]
and the $\P$-functors
\[ F_k\colon D^b(B) \to D^b(X,\beta) \]
are well-defined, and the $\P$-twist associated to $F_k$ satisfies \eqref{flop_flop_fmla}.
\end{prop}

\begin{proof}
Let $\bar X$ be the space obtained from $X$ by contracting $P$ down to $B$, or from $X'$ by contracting $P^*$ down to $B$:
\[ \xymatrix{
X \ar[rrdd]_f & P \ar@{_(->}[l]_j \ar[rd]_\varpi
& & P^* \ar@{^(->}[r]^{j'} \ar[ld]^{\varpi'} & X' \ar[lldd]^{f'} \\
& & B \ar@{^(->}[d]^{\bar\jmath} \\
& & {\phantom.}\bar X.
} \]
Note that $\hat X = X \times_{\bar X} X'$.  We work in the analytic category to avoid worrying about whether $\bar X$ is projective.

We claim that
\begin{equation} \label{br_pullback}
f^*\colon \Br(\bar X) \to \Br(X)
\end{equation}
is an isomorphism.  The proof is identical to the one for a blow-up along a smooth center.\footnote{Cf.\ Grothendieck's classic proof \cite[Thm.~7.1]{brauer3} in the \'etale topology.}  Taking the exponential sequence
\[ 0 \to \Z \to \O_X \to \O_X^* \to 0 \]
and pushing down to $\bar X$, we find that
\begin{align*}
R^0 f_*(\O_X^*) &= \O_{\bar X}^* \\
R^1 f_*(\O_X^*) &= \bar\jmath_*(\Z) \\
R^2 f_*(\O_X^*) &= 0.
\end{align*}
Then from the Leray spectral sequence we get an exact sequence
\[ 0 \to \Pic(\bar X) \xrightarrow{f^*} \Pic(X) \to \Z \to H^2(\O_{\bar X}^*) \xrightarrow{f^*} H^2(\O_X^*) \to H^1(B,\Z), \]
where the map $\Pic(X) \to \Z$ takes a line bundle on $X$, restricts it to $P$, and asks for its degree on a fiber of $\varpi$.  By hypothesis this map is surjective, so the map $f^*\colon H^2(\O_{\bar X}^*) \to H^2(\O_X^*)$ is injective.  Taking torsion parts and noting that $H^1(B,\Z)$ is torsion-free, we see that \eqref{br_pullback} is an isomorphism.

Thus we get a class $\bar\beta \in \Br(\bar X)$ such that $\beta = f^* \bar\beta$ and $\beta' = f'^* \bar\beta$.  We choose once and for all a cocycle representing $\bar\beta$; this determines cocycles representing $\beta$ and $\beta'$ such that the cocycle representing $(\beta^{-1} \boxtimes \beta')|_{\hat X}$ is canonically trivial.  Thus the pushforward map
\[ D^b(\hat X) \to D^b(X \times X', \beta^{-1} \boxtimes \beta') \]
is canonically defined, and we can take $\KN_k$ to be the functor $D^b(X,\beta) \to D^b(X',\beta')$ induced by the pushforward of $\L^{\otimes k}$.

Next we turn to the $\P$-functor and its $\P$-twist.  We have the natural inclusion $\varpi \times j\colon P \to B \times X$, and the cocycle representing $(\bar\jmath^* \bar\beta^{-1} \boxtimes \beta)|_P$ is again canonically trivial, so the pushforward of $\O_{P/B}(k)$ induces a $\P$-functor
\[ F_k\colon D^b(B, \bar\jmath^* \bar\beta) \to D^b(X,\beta) \]
satisfying the same formula \eqref{flop_flop_fmla} as in the untwisted case.  It remains to compare $D^b(B, \bar\jmath^* \bar\beta)$ to $D^b(B)$.  By hypothesis we have
\[ \varpi^* \bar\jmath^* \bar\beta = j^* f^* \bar\beta = j^* \beta = 1. \]
The pullback $\varpi^*\colon \Br(B) \to \Br(P)$ is an isomorphism because $P$ is the projectivization of a vector bundle; thus $\bar\jmath^* \bar\beta$ is trivial.  But we have not arranged for it to be canonically trivial: that is, the 2-cocycle representing $\bar\jmath^* \bar\beta$ is only the coboundary of a 1-cocycle, which we need to choose.  Any choice gives an equivalence $D^b(B) \to D^b(B, \bar\jmath^* \bar\beta)$, and thus a $\P$-functor $F_k\colon D^b(B) \to D^b(X,\beta)$.  If we make a different choice then the equivalence may differ by tensoring by a line bundle on $B$.  But by fact \opt{ams}{\eqref{third_fact}}\opt{lms}{\ref{third_fact}} from earlier this gives an isomorphic twist $P_k$, so \eqref{flop_flop_fmla} continues to hold.
\end{proof}

In the proof of Proposition \ref{fma_thm}\opt{ams}{(b)}\opt{lms}{(ii)}, we will see that our ``mild hypothesis'' holds in our application, where $P = \C$ and $X = \Picbar^{-1}$.

% !TeX root = main.tex

\section{Factorization} \label{fma}

Recall from the introduction that $S$ is a K3 surface and $\O_S(1)$ is an ample generator of $\Pic(S)$, of degree $2g-2$.  Thus $V := H^0(\O_S(1))$ is $(g+1)$-dimensional,\footnote{For this and other general facts about curves on K3 surfaces we recommend \cite[Ch.~2]{huybrechts_K3}.} and the linear system $|\O_S(1)| = \P V \cong \P^g$.  The general member of $\P V$ is a smooth curve of genus $g$; moreover, since $\O_S(1)$ generates $\Pic(S)$, every member of the linear system is reduced and irreducible, so we can use Arinkin's results \cite{arinkin} on compactified Jacobians of integral curves with planar singularities.  We introduce some more notation:
\smallskip
\setlength \leftmargini {1.5em}
\begin{itemize}
\setlength \itemsep {3pt}
\item $\C = \{ (x,C) \in S \times \P V : x \in C \}$ is the universal curve.
\item $\varpi\colon \C \to S$ is the projection $\varpi(x,C) = x$, which is a $\P^{g-1}$-bundle.
\item $\Picbar^d = \Picbar^d(\C/\P V)$ is the relative compactified Jacobian, or equivalently the moduli space of stable sheaves on $S$ with Mukai vector $(0,1,d+1-g)$.
\item $\alpha_d \in \Br(\Picbar^d)$ is the Brauer class obstructing the existence of a universal sheaf on $S \times \Picbar^d$.  Because the Hilbert polynomial of the sheaves is
\[ (2g-2)t + (d+1-g) \]
we have $\alpha_d^{2g-2} = \alpha_d^{d+1-g} = 1$.  In particular, if $2g-2$ is relatively prime to $d+1-g$ then $\alpha_d = 1$ and $\Picbar^d$ is a fine moduli space.
\item $\F_d$ is a $(1 \boxtimes \alpha_d)$-twisted pseudo-universal sheaf on $S \times \Picbar^d$.  Observe that this is supported on $\C \times_{\P V} \Picbar^d$, where the latter is embedded in $S \times \Picbar^d$ via $\varpi \times 1$.\footnote{The reader may remark that $\F_d$ has rank 1 on $\C \times_{\P V} \Picbar^d$, so the Brauer class $(1 \boxtimes \alpha_d)$ on $S \times \Picbar^d$ must become trivial when restricted $\C \times_{\P V} \Picbar^d$.  Why then do we not ``untwist'' $\F_d$ to get an honest sheaf?  Because it would no longer be a universal sheaf: we would have to pass to an open cover of $\C \times_{\P V} \Picbar^d$ whose open sets were not unions of fibers of $\C \to \P V$, and the re-glued $\F_d$ would no longer parametrize the correct sheaves on those fibers.  \Caldararu\ has discussed the same issue for $g=0$ in \cite[\S4.3]{andrei}.  Lieblich has inveighed against untwisting in \cite[Rmk.~1.3.3]{lieblich}.}
\item $AJ\colon \C \hookrightarrow \Picbar^{-1}$ is the Abel--Jacobi embedding, sending $(x,C)$ to the ideal sheaf $\I_{x/C}$; this is a stable sheaf because $C$ is reduced and irreducible.  In the proof of Proposition \ref{fma_thm}\opt{ams}{(b)}\opt{lms}{(ii)} we will see that $AJ^*(\alpha_{-1}) = 1$.
% note that \alpha_{-1} is already trivial if g is odd.
\end{itemize}
\setlength \leftmargini {2em}
\bigskip % This and the next \bigskip are mainly to get the proof of Prop 2.1 onto the next page.

\begin{prop} \label{fma_thm} With the notation set up in the previous paragraph:
\begin{enumerate}
\item For every $m$ and $n$ there is an $(\alpha_m^{-n} \boxtimes \alpha_n^{-m})$-twisted Poincar\'e sheaf $\bar P_{mn}$ on $\Picbar^m \times_{\P V} \Picbar^n$, such that the induced map
\[ F\!M \colon D^b(\Picbar^m, \alpha_m^n) \to D^b(\Picbar^n, \alpha_n^{-m}) \]
is an equivalence.
\item Consider the embedding
\[ AJ \times 1\colon \, \C \times_{\P V} \Picbar^n \, \hookrightarrow \, \Picbar^{-1} \times_{\P V} \Picbar^n. \]
Given any trivialization of $AJ^*(\alpha_{-1})$, we can construct $\bar P_{-1n}$ so that
\[ (AJ \times 1)^* \bar P_{-1n} \cong \F_n \otimes u_1^* \varpi^* \L, \]
where $u_1$ is the projection of $\C \times_{\P V} \Picbar^n$ onto the first factor and $\L$ is a line bundle on $S$.  This implies the factorization \eqref{factor_F} from the introduction.
\end{enumerate}
\end{prop}
\bigskip % This and the previous \bigskip are mainly to get the proof of Prop 2.1 onto the next page.
\opt{lms}{\begin{proof*}}\opt{ams}{\begin{proof}}
\opt{ams}{(a)}\opt{lms}{(i)} First we recall Arinkin's description of the Poincar\'e line bundle for a single curve \cite[eq.~(1.1)]{arinkin}.  At a point
\[ (F,G) \in (\Pic^0(C) \times \Picbar^0(C)) \cup (\Picbar^0(C) \times \Pic^0(C)), \]
the fiber is\footnote{In fact this is the dual of Arinkin's definition, but this does not cause any problems; see \cite[Lem.~6.2]{arinkin}.}
\begin{equation} \label{arinkin_poincare}
\det \RGamma(F \otimes G)^{-1} \otimes \det \RGamma(\O_C)^{-1} \otimes \det \RGamma(F) \otimes \det \RGamma(G).
\end{equation}
Note that $F \otimes G$ is a sheaf since at least one of $F$ and $G$ is a line bundle; if both were sheaves then $F \otimes G$ might be an unbounded complex, so the first term of \eqref{arinkin_poincare} would not be well-defined.

Next we globalize this to our family of curves.  For the reader's convenience we display the big diagram:
\[ \xymatrix{
& \C \times_{\P V} \Picbar^m \times_{\P V} \Picbar^n \ar[ld]_{p_{12}} \ar[d]_{p_{13}} \ar[rd]^{p_{23}} \\
\C \times_{\P V} \Picbar^m \ar[d]_{q_1} \ar[rd]_<{q_2}
& \C \times_{\P V} \Picbar^n \ar[ld]^<{r_1}|\hole \ar[rd]_<{r_2}|\hole
& \Picbar^m \times_{\P V} \Picbar^n \ar[ld]^<{u_1} \ar[d]^{u_2} \\
\C \ar[rd]_{v_1} & \Picbar^m \ar[d]_{v_2} & \Picbar^n \ar[ld]^{v_3} \\
& \P V
} \]
Following \eqref{arinkin_poincare} we define a twisted line bundle $P_{mn}$ on
\begin{equation} \label{Pic_x_Pic_natural}
(\Pic^m \times_{\P V} \Picbar^n) \cup (\Picbar^m \times_{\P V} \Pic^n)
\end{equation}
by the formula
\begin{multline*}
P_{mn} = (\det p_{23*} (p_{12}^* \F_m \otimes p_{13}^* \F_n))^{-1}
\otimes u_1^* v_2^* (\det v_{1*} \O_\C)^{-1} \\
\otimes u_1^*(\det q_{2*} \F_m) \otimes u_2^*(\det r_{2*} \F_n).
\end{multline*}
We will see that this is $(\alpha_m^{-n} \boxtimes \alpha_n^{-m})$-twisted.  Then we will define $\bar P_{mn} = \iota_* P_{mn}$, where $\iota$ is the inclusion of \eqref{Pic_x_Pic_natural} in $\Picbar^m \times_{\P V} \Picbar^n$.  Part (a) asserts that this is coherent, and that its convolution with its left and right adjoint kernels is $\O_\Delta$.  Arinkin has proved this for $m=n=0$; to extend it to all $m$ and $n$ we will argue that locally in $\P V$, our $\bar P_{mn}$ differs from $\bar P_{00}$ by a line bundle of the form $\L' \boxtimes \L''$.

Now we fill in the details of the outline just given.  Begin by choosing an analytic (or \'etale) open cover $\{ U_i \}$ of $\P V$ over which there are local sections $s_i$ of the universal curve $\C \to \P V$.  Regard $s_i$ as a divisor in $\C|_{U_i}$, and observe that it necessarily stays in the smooth locus of each fiber, so $\O_{\C|_{U_i}}\!(s_i)$ is a line bundle.

Because $\C|_{U_i}$ has a section, there is a universal sheaf $\F_{m,i}$ on $(\C \times_{\P V} \Picbar^m)|_{U_i}$, and similarly with $n$.  One way to see this is as follows.  Fix a $d$ such that $\Picbar^d$ is a fine moduli space, so there is an honest universal sheaf $\F_d$ on $\C \times_{\P V} \Picbar^d$.  Let
\[ t_{(d-m) s_i} \colon \Picbar^m|_{U_i} \to \Picbar^d|_{U_i} \]
be the ``translation'' isomorphism, which sends a sheaf $F$ on a curve $C$ to $F \otimes \O_C((d-m) s_i)$.  Then set
\begin{equation} \label{def_of_F}
\F_{m,i} = (1 \times t_{(d-m) s_i})^* \F_d \otimes q_1^* \O_{\C|_{U_i}}\!((m-d)s_i).
\end{equation}
This is one possible construction of $\F_{m,i}$; any other construction differs from this by a line bundle pulled back from $\Picbar^m|_{U_i}$.

Let $U_{ij} = U_i \cap U_j$ as usual.  Then $\F_{m,i}|_{U_{ij}}$ and $\F_{m,j}|_{U_{ij}}$ are two universal sheaves, so they differ by a line bundle $\L_{m,ij}$ pulled back from $\Picbar^m|_{ij}$:
\begin{equation} \label{F_transforms}
\F_{m,j}|_{U_{ij}} = \F_{m,i}|_{U_{ij}} \otimes q_2^* \L_{m,ij}.
\end{equation}
If we define $\F_{m,i}$ as in \eqref{def_of_F} then we can give an explicit formula for $\L_{m,ij}$ involving $s_j - s_i$, but we will not need it.  The line bundles $\L_{m,ij}$ represent the Brauer class $\alpha_m$ as a gerbe \`a la Hitchin; see \cite[p.~13]{andrei}.\footnote{To get an $\alpha_m$-twisted sheaf in the more familiar sense of \cite[Def.~1.2.1]{andrei} we would trivialize $\L_{m,ij}$ over a finer open cover of $\Picbar^m$ and go from there, but this would obscure the picture needlessly.}

Now we define $P_{mn,i}$ on the restriction of \eqref{Pic_x_Pic_natural} to $U_i$ by the same formula as above:
\begin{multline*}
P_{mn,i} = (\det p_{23*} (p_{12}^* \F_{m,i} \otimes p_{13}^* \F_{n,i}))^{-1}
\otimes u_1^* v_2^* (\det v_{1*} \O_{\C|_{U_i}})^{-1} \\
\otimes u_1^*(\det q_{2*} \F_{m,i}) \otimes u_2^*(\det r_{2*} \F_{n,i}).
\end{multline*}
We study how this transforms under various changes to $\F_{m,i}$ and $\F_{n,i}$.  First, if we tensor $\F_{m,i}$ by a line bundle pulled back from $\Picbar^m$ then we find that
\begin{align} \label{lb_for_m}
\F_{m,i} &\rightsquigarrow \F_{m,i} \otimes q_2^* \L &\Longrightarrow& & P_{mn,i} &\rightsquigarrow P_{mn,i} \otimes u_1^* \L^{-n},
\end{align}
and similarly
\begin{align} \label{lb_for_n}
\F_{n,i} &\rightsquigarrow \F_{n,i} \otimes r_2^* \L &\Longrightarrow& & P_{mn,i} &\rightsquigarrow P_{mn,i} \otimes u_2^* \L^{-m},
\end{align}
in agreement with \cite[eq.~(4.2)]{mrv2}.  Second, suppose that we tensor $\F_{m,i}$ with $q_1^* \O_{\C|_{U_i}}\!(-s_i)$.  Let $\F_{m,i}|_{s_i}$ denote the line bundle on $\Picbar^m|_{U_i}$ which might more properly be called $q_{2*}(\F_{m,i} \otimes q_1^* \O_{s_i})$.\footnote{Or $(s_i \times 1)^* \F_{m,i}$, if we were regarding $s_i$ as a map rather than a divisor.}  Then using the exact sequence
\[ 0 \to \O_{\C|_{U_i}}\!(-s_i) \to \O_{\C|_{U_i}} \to \O_{s_i} \to 0 \]
we find that
\begin{align} \label{deg_change_for_m}
\F_{m,i} &\rightsquigarrow \F_{m,i} \otimes q_1^* \O_{\C|_{U_i}}\!(-s_i) &\Longrightarrow& & P_{mn,i} &\rightsquigarrow P_{mn,i} \otimes u_2^* \F_{n,i}|_{s_i},
\end{align}
and similarly
\begin{align} \label{deg_change_for_n}
\F_{n,i} &\rightsquigarrow \F_{n,i} \otimes r_1^* \O_{\C|_{U_i}}\!(-s_i) &\Longrightarrow& & P_{mn,i} &\rightsquigarrow P_{mn,i} \otimes u_1^* \F_{m,i}|_{s_i}.
\end{align}

Now from \eqref{F_transforms}, \eqref{lb_for_m}, and \eqref{lb_for_n} we see that
\[ P_{mn,j}|_{U_{ij}} = P_{mn,i}|_{U_{ij}} \otimes (\L_{m,ij}^{-n} \boxtimes \L_{n,ij}^{-m}), \]
so the line bundles $P_{mn,i}$ glue to give an $(\alpha_m^{-n} \boxtimes \alpha_n^{-m})$-twisted line bundle $P_{mn}$ on \eqref{Pic_x_Pic_natural}, as claimed.

It remains to show that $\bar P_{mn,i}$ differs from $\bar P_{00,i}$, or more properly from $(t_{-m s_i} \times t_{-n s_i})^* \bar P_{00,i}$, by a line bundle of the form $\L' \boxtimes \L''$.  From \eqref{def_of_F} and the line after it we see that $\F_{m,i}$ differs from $(1 \times t_{-m s_i})^* \F_0$ by a line bundle of the form $\O_{\C|_{U_i}}\!(m s_i) \boxtimes \L$.  Thus from \eqref{lb_for_m}, \eqref{lb_for_n}, \eqref{deg_change_for_m}, and \eqref{deg_change_for_n} we see that $P_{mn,i}$ differs from $(t_{-m s_i} \times t_{-n s_i})^* P_{00,i}$ by a line bundle of the form $\L' \boxtimes \L''$.  By the projection formula, the same holds for $\bar P_{mn,i}$ and $(t_{-m s_i} \times t_{-n s_i})^* \bar P_{00,i}$.
\medskip

\opt{ams}{(b)}\opt{lms}{(ii)} Repeat the big diagram from above but with $\Picbar^m$ replaced by $\C$:
\[ \xymatrix{
& \C \times_{\P V} \C \times_{\P V} \Picbar^n \ar[ld]_{p_{12}} \ar[d]_{p_{13}} \ar[rd]^{p_{23}} \\
\C \times_{\P V} \C \ar[d]_{q_1} \ar[rd]_<{q_2}
& \C \times_{\P V} \Picbar^n \ar[ld]^<{r_1}|\hole \ar[rd]_<{r_2}|\hole
& \C \times_{\P V} \Picbar^n \ar[ld]^<{u_1} \ar[d]^{u_2} \\
\C \ar[rd]_{v_1} & \C \ar[d]_{v_2} & \Picbar^n \ar[ld]^{v_3} \\
& \P V
} \]
This maps to the earlier diagram with $m=-1$ using the Abel--Jacobi map
\[ AJ\colon \C \hookrightarrow \Picbar^{-1}. \]
The Abel--Jacobi map is the classifying map for the ideal sheaf of the diagonal in $\C \times_{\P V} \C$.  In particular $AJ^*(\alpha_{-1})$ is trivial, and after choosing a trivialization of it we can pull back $\alpha_{-1}$-twisted sheaves on $\Picbar^{-1}$ to untwisted sheaves on $\C$; let such a trivialization be given.  Then there is a line bundle $\N$ on $\C$ such that
\[ (1 \times AJ)^* \F_{-1} \cong \I_\Delta \otimes q_2^* \N \]
on $\C \times_{\P V} \C$.  Recall that we have
\[ \xymatrix{
\opt{lms}{\,}\C\opt{lms}{\,} \ar[d]_{\varpi} \ar[rd]^{v_2} \ar@{^(->}[rr]^-{AJ} & & \Picbar^{-1} \ar[ld]_{v_2} \\
S & \:\P V,
} \]
where $\varpi$ is a $\P^{g-1}$-bundle and $v_2\colon \C \to \P V$ embeds the fibers of $\varpi$ as hyperplanes in $\P V = \P^g$.  Thus we can write
\[ \N = \varpi^* \N' \otimes v_2^* \O_{\P V}(k) \]
for some line bundle $\N'$ on $S$ and some $k \in \Z$, and we can absorb the $\O_{\P V}(k)$ into $\F_{-1}$, leaving us with
\begin{equation} \label{L_prime_eq}
(1 \times AJ)^* \F_{-1} \cong \I_\Delta \otimes q_2^* \varpi^* \N'.
\end{equation}
% With a lot more work we could absorb the L' as well, but it doesn't seem worth it.

Suppose for a moment that $\N'$ is trivial, so we have
\[ (1 \times AJ)^* \F_{-1} \cong \I_\Delta. \]
By a series of base changes we find that on the rightmost $\C \times_{\P V} \Picbar^n$,
\begin{multline} \label{delta_thing}
(AJ \times 1)^* P_{-1n} \cong (\det p_{23*} (p_{12}^* \I_\Delta \otimes p_{13}^* \F_n))^{-1}
\otimes u_1^* v_2^* (\det v_{1*} \O_\C)^{-1} \\
\otimes u_1^*(\det q_{2*} \I_\Delta) \otimes u_2^*(\det r_{2*} \F_n).
\end{multline}
Note that \eqref{delta_thing} is only defined on
\begin{equation} \label{C_x_Pic_natural}
(\C^\circ \times_{\P V} \Picbar^n) \cup (\C \times_{\P V} \Pic^n),
\end{equation}
where $\C^\circ \subset \C$ is the smooth locus of the map $v_1\colon \C \to \P V$, i.e.\ the union of the smooth loci of all the curves.  Using the exact sequence
\[ 0 \to \I_\Delta \to \O \to \O_\Delta \to 0 \]
we can simplify \eqref{delta_thing}: after a long exercise in base change and the projection formula, the first term becomes $u_2^*(\det r_{2*} \F_n)^{-1} \otimes \det \F_n$, which partially cancels with the fourth term to leave $\det \F_n$, and the third term becomes $u_1^* v_2^* (\det v_{1*} \O_\C)$ which cancels with the second term; thus
\[ (AJ \times 1)^* P_{-1n} \cong \det \F_n. \]
Since $\F_n$ is a (twisted) line bundle on \eqref{C_x_Pic_natural}, this equals $\F_n$.

Now as in the proof of \cite[Lem.~6.4]{arinkin} we observe that $(AJ \times 1)^* \bar P_{-1n}$ and $\F_n$ are both maximal Cohen--Macaulay sheaves whose restrictions to \eqref{C_x_Pic_natural} agree, so they agree.  For a review of the necessary Cohen--Macaulay machinery we recommend \cite[\S2]{arinkin}.

Finally, if the line bundle $\N'$ of \eqref{L_prime_eq} is non-trivial then as in \eqref{lb_for_m} we find that
\[ \opt{lms}{\singlebox} (AJ \times 1)^* \bar P_{-1n} \cong \F_n \otimes u_1^* \varpi^* \N'^{-n}. \opt{lms}{\esinglebox} \opt{ams}{\qedhere} \]
\opt{lms}{\end{proof*}}\opt{ams}{\end{proof}}

\begin{rmk}
The proof above is valid for any surface, and for any base-point-free linear system of curves (of any genus) whose members are all reduced and irreducible.  It should be straightforward to extend part \opt{ams}{(a)}\opt{lms}{(i)} to linear systems with reducible (but still reduced) members, and in particular to primitive divisor classes on K3 surfaces of higher Picard rank, using \cite{mrv2}.  But then part \opt{ams}{(b)}\opt{lms}{(ii)}, which is essential to our application, falls apart: the Abel--Jacobi map may fail to exist, as the ideal sheaf of a point in a reducible curve may fail to be stable.
\end{rmk}

\begin{rmk} \label{special_cases}
We point out some interesting special cases of Proposition \ref{fma_thm}.   First, for all $d$ we have
\[ D^b(\Picbar^d) \cong D^b(\Picbar^0, \alpha_0^d). \]
Second, we have
\begin{align*}
D^b(\Picbar^0) &\cong D^b(\Picbar^{g-1}), \\
D^b(\Picbar^1) &\cong D^b(\Picbar^{g-2}), \\
D^b(\Picbar^2) &\cong D^b(\Picbar^{g-3}),
\end{align*}
and so on, because $\alpha_d^{g-1-d} = 1$ and $\alpha_{g-1-d}^d = 1$.
\end{rmk}

\begin{rmk}
The \emph{period} of a compact hyperk\"ahler manifold $X$ is the second cohomology group $H^2(X,\Z)$ with its weight-2 Hodge structure and Beauville--Bogomolov pairing.  It is a birational invariant of compact hyperk\"ahler manifolds \cite[Cor.~4.7]{huybrechts_basic}, but as a consequence of the previous remark we note that it is not a derived invariant:

{\renewcommand \thethm {\ref{per_not_inv}}
\addtocounter{thm}{-1}
\begin{thm}
For every $g \ge 2$ there are compact hyperk\"ahler $2g$-folds $X$ and $Y$ such that $D^b(X) \cong D^b(Y)$ but $H^2(X,\Z)$ is not Hodge-isometric to $H^2(Y,\Z)$.  In particular $X$ is not birational to $Y$.
\end{thm}}
\begin{proof}
Take $X = \Picbar^0$ and $Y = \Picbar^{g-1}$.  These are derived equivalent by Remark \ref{special_cases}.   Following Sawon \cite[Proof of Prop.\ 15]{sawon} % Prop. 13 in the arXiv version.
we find that the Picard lattice of $X$ has discriminant $-4$, whereas the Picard lattice of $Y$ has discriminant $-1$.
\end{proof}

The corresponding statement for K3 surfaces is well-known: there are many examples of K3 surfaces that are derived equivalent but not birational (i.e.\ not isomorphic).  For Calabi--Yau 3-folds, Borisov and \Caldararu\ \cite{bc} produced the first example of a pair that are derived equivalent but not birational; but derived-equivalent Calabi--Yau 3-folds necessarily have Hodge-isometric periods $H^3(-,\Z)/$torsion.\footnote{This can be seen as follows: if $X$ is a Calabi--Yau 3-fold then its Atiyah--Hirzebruch spectral sequence degenerates at the $E_2$ page, so $H^3(X,\Z)/$torsion is the image of the Mukai vector $v\colon K^1_\text{top}(X) \to H^3(X,\mathbb Q)$, and this is a derived invariant.  For more detail and references see \cite[Proof of Prop.]{br_not_inv} and \cite[\S2]{at}.  Note that \Caldararu\ \cite[Prop.~3.1]{andrei_7sec} showed by an easier argument that $H^3(X, \Z[\tfrac12])/$torsion is a derived invariant.}

It is not quite true that if two hyperk\"ahler varieties have Hodge-isometric periods then they are birational; for varieties deformation-equivalent to moduli spaces of sheaves on K3 surfaces, there is a slightly larger lattice $\tilde\Lambda \supset H^2(X,\Z)$, the \emph{Markman--Mukai lattice} \cite[\S9]{survey_of_torelli}, which controls the birational geometry.  It is interesting to note that in our example, the Markman--Mukai lattices of $X$ and $Y$ are isomorphic: $\tilde\Lambda \cong H^*(S,\Z)$ in both cases.  Two K3 surfaces are derived equivalent if and only if they have the same Mukai lattice \cite{orlov}; in higher dimensions it is not clear whether one should expect the Markman--Mukai lattice to be a derived invariant.
\end{rmk}

\begin{rmk}
As D.~Huybrechts has pointed out to us, it is very likely that for
every $n$ there are K3 surfaces $S$ and $S'$ such that $D^b(S) \cong D^b(S')$, so $D^b(\Hilb^n(S)) \cong D^b(\Hilb^n(S'))$ by \cite[Prop.~8]{ploog}, but $\Hilb^n(S)$ and $\Hilb^n(S')$ are not birational, giving a simpler example of Theorem \ref{per_not_inv}.  Thanks to the Torelli theorems of Mukai--Orlov and Verbitsky--Markman, this is a purely lattice-theoretic question, and it should be answerable using the methods of \cite{hloy,stellari,hp}, but to our knowledge no such analysis has appeared in the literature.
\end{rmk}

\begin{rmk}
It follows from Arinkin's autoduality theorem \cite[Thm.~B]{arinkin} that the moduli space $\Picbar^0({\Picbar^n}/\P V)$ is isomorphic to $\Picbar^0 = \Picbar^0(\C/\P V)$ for any $n$.  Inspecting the formula for the $(1 \boxtimes \alpha_0^{-n})$-twisted sheaf $\bar P_{n0}$ on $\Picbar^n \times_{\P V} \Picbar^0$, we see that it has degree 0 on the fibers of the projection to $\Picbar^0$, and indeed that it is universal for this moduli problem; thus the Brauer class for this moduli problem is $\alpha_0^{-n}$.

Since $\alpha_0^{1-g} = 1$ but $\Picbar^{g-1} \not\cong \Picbar^0$, this disagrees with a proposition of Sawon \cite[Prop.~9]{sawon}. % Prop. 17 in the arXiv version, but the published version is more detailed.
He considers abelian fibrations $X \to B$ satisfying certain hypotheses, and lets $P = \Picbar^0(X/B)$ and $X^0 = \Picbar^0(P/B)$.  In our example we have $B = \P V$, $X = \Picbar^{g-1}(\C/\P V)$, and $P = X^0 = \Picbar^0(\C/\P V)$.  He states that the following are equivalent: (1) $X$ is isomorphic to $X^0$ over $B$, (2) $X \to B$ admits a section, (3) there is a global universal sheaf on $X \times P$, and (4) a certain Brauer class $\beta$ on $P$ vanishes.  In fact one has $(1) \Leftrightarrow (2) \Rightarrow (3) \Leftrightarrow (4)$, but not $(4) \Rightarrow (1)$.
% Had more detail on this but moved it to a separate file.
\end{rmk}

% !TeX root = main.tex

\section{Compatibility with the Hilbert scheme} \label{compatibility}

We continue the notation of the previous section: thus $S$ is a K3 surface with $\Pic(S) = \Z$ generated by an ample line bundle $\O_S(1)$ of degree $2g-2$, and $V = H^0(\O_S(1)) = \CC^{g+1}$, so we get a natural map $S \to \P V^*$ and a universal curve $\C \to \P V$ whose general fiber is a smooth curve of genus $g$.

Consider the moduli spaces
\begin{align*}
\Hilb^g &:= \Hilb^g(S) = \M(1,0,1-g) \\
\Picbar^{-g} &:= \Picbar^{-g}(\C/\P V) = \M(0,1,1-2g)
\end{align*}
and the birational map $\Hilb^g \dashrightarrow \Picbar^{-g}$ discussed in the introduction: for a generic length-$g$ subscheme $\zeta \subset S$, the image of $\zeta$ in $\P V^* = \P^g$ spans a hyperplane, so $\zeta$ is contained in a unique curve $C \in \P V$, and we map $\zeta \in \Hilb^g$ to $\O_C(-\zeta) := \I_{\zeta/C} \in \Picbar^{-g}$.  This birational map is resolved by the correspondence
\[ \tilde X := \Hilb^g(\C/\P V) = \{ (\zeta, C) : \zeta \subset C \} \subset {\Hilb^g} \times \P V, \]
in the following way:
\[ \xymatrix{
& \tilde X \ar[ld]_p \ar[rd]^q & & & & \zeta \subset C \ar@{|->}[ld] \ar@{|->}[rd] \\
\Hilb^g & & \Picbar^{-g} & & \I_{\zeta/S} & & {\phantom.}\I_{\zeta/C}.
} \]
Let $P^* \subset \Hilb^g$ be the indeterminacy locus of $p^{-1}$, consisting of subschemes $\zeta \subset S$ for which $h^0(\I_{\zeta}(1))$ jumps from 1 to 2 or more.  Let $P \subset \Picbar^{-g}$ be the indeterminacy locus of $q^{-1}$, consisting of sheaves $\ell$ on curves $C$ for which $\operatorname{hom}(\ell, \O_C)$ jumps from 1 to 2 or more.  Both of these have codimension 2.  If $g \le 5$ then by \cite[Example~21]{brill_noether}, the jumping is to 2 and not more, the loci $P$ and $P^*$ are dual $\P^2$-bundles over a smooth hyperk\"ahler $(2g-4)$-fold $B := \M(2,1,1)$,
% To use this reference notice that \M(1,0,1-g) = \M(1,1,0).
and $q \circ p^{-1}$ is a Mukai flop.

As we said in the introduction, $q \circ p^{-1}$ is implemented by the spherical twist around $\O_S(-1)$, in the following sense: if $\zeta \notin P^*$, so $\zeta$ is contained in a unique curve $C$, then we have
\begin{align}
T_{\O_S(-1)}(\I_\zeta)
&= \cone\!\big(\; \O_S(-1) \otimes \RGamma(\I_\zeta(1)) \to \I_\zeta \;\big) \label{line_with_T_O(-1)} \\
&= \cone\!\big(\; \I_{C/S} \to \I_{\zeta/S} \;\big) \notag \\
&= \I_{\zeta/C}. \notag
\end{align}
Letting $\zeta$ vary, this suggests that on $S \times \Hilb^g$ there should be some relation between
\begin{equation} \label{cone_we're_taking}
\cone\!\big(\; \O_S(-1) \boxtimes \pi_{\Hilb,*}(\I_Z \otimes \pi_{S}^* \O_S(1)) \longrightarrow \I_Z \;\big),
\end{equation}
where $Z \subset S \times \Hilb^g$ is the universal subscheme, and
\[ p_* q^* \F \]
where $\F$ is a universal sheaf on $S \times \Picbar^{-g}$.\footnote{More properly we should write $(1 \times p)_* (1 \times q)^* \F$, but such heavy notation would become untenable in what follows.}  In Proposition \ref{compat_prop} below we make this precise.  In Proposition \ref{intertwinement_prop} we use it to compare the $\P^{g-1}$-functors $D^b(S) \to D^b(\Hilb^g)$ and $D^b(S) \to D^b(\Picbar^{-g})$ induced by $\I_Z$ and $\F$, and the associated $\P$-twists.

\begin{prop} \label{compat_prop}
For a suitable normalization of the universal sheaf $\F$ on $S \times \Picbar^{-g}$, the cone \eqref{cone_we're_taking} is isomorphic to $p_* (q^* \F(2E)),$ where $E \subset \tilde X$ is the exceptional divisor of the blow-up $q\colon \tilde X \to \Hilb^g$ or $p\colon \tilde X \to \Picbar^{-g}$.
\end{prop}
\begin{proof}
First we study
\begin{equation} \label{underlined_term} % In an earlier version it was underlined.
\pi_{\Hilb,*}(\I_Z \otimes \pi_S^* \O_S(1))
\end{equation}
which appears in \eqref{cone_we're_taking}.  Let $W$ denote the rank-$g$ vector bundle on $\Hilb^g$ whose fiber over a point $\zeta$ is $H^0(\O_\zeta(1))$.  Take the exact sequence
\[ 0 \to \I_Z \to \O_{S \times \Hilb} \to \O_Z \to 0, \]
tensor with $\pi_S^* \O_S(1)$, and push down to $\Hilb^g$ to get an exact triangle
\begin{equation} \label{underlined_term_seq}
\pi_{\Hilb,*}(\I_Z \otimes \pi_S^* \O_S(1)) \to \O_{\Hilb} \otimes V \to W.
\end{equation}
The map $\O_{\Hilb} \otimes V \to W$ is surjective away from $P^* \subset \Hilb^g$, and for later use we record the following lemma.
\begin{lem}
The variety $\tilde X \subset {\Hilb^g} \times \P V$ is cut out by the transverse section of $W \boxtimes \O_{\P V}(1)$ corresponding to the map $\O_{\Hilb} \otimes V \to W$ of \eqref{underlined_term_seq}.
\end{lem}
\begin{proof}
Observe that $\tilde X$ is the locus where the composition
\[ \O_{\Hilb} \boxtimes \O_{\P V}(-1) \to \O_{{\Hilb} \times \P V} \otimes V \to W \boxtimes \O_{\P V} \]
vanishes, and this composition amounts to a section of $W \boxtimes \O_{\P V}(1)$.  Since $\operatorname{codim} \tilde X = \rank W$, the section is transverse.
\end{proof}

\noindent \emph{Proof of Proposition \ref{compat_prop}, continued.}
At a point $(\zeta, C) \in \tilde X$, we have an exact sequence of sheaves on $S$
\[ 0 \to \I_{C/S} \to \I_{\zeta/S} \to \I_{\zeta/C} \to 0, \]
which we can rewrite as
\[ 0 \to \O_S(-1) \to \I_\zeta \to \O_C(-\zeta) \to 0. \]
We wish to write the family version of this on $S \times \tilde X$.  We have
\[ S \times \tilde X = \{ (x,\zeta,C) : \zeta \subset C \} \subset S \times {\Hilb^g} \times \P V. \]
Inside this consider
\begin{align*}
Z' &:= \{ (x,\zeta,C) : x \in \zeta \subset C) \} \\
\C' &:= \{ (x,\zeta,C) : x \in C, \zeta \subset C \}.
\end{align*}
We have $Z' \subset \C' \subset S \times \tilde X$, so we get an exact sequence
\[ 0 \to \I_{\C'/S \times \tilde X} \to \I_{Z'/S \times \tilde X} \to \I_{Z'/\C'} \to 0. \]
We claim this can be rewritten as
\begin{equation} \label{main_seq}
0 \to \O_S(-1) \boxtimes r^* \O_{\P V}(-1) \to p^* \I_Z \to q^* \F (E) \to 0,
\end{equation}
where $r\colon \tilde X \to \P V$ is the natural map, for a suitable normalization of $\F$.
\begin{enumerate}
\item For the first term, observe that the divisor $\C' \subset S \times \tilde X$ is the pullback of the divisor $\C \subset S \times \P V$, whose associated line bundle is $\O_{S \times \P V}(1,1)$.
\item For the second term, observe that the intersection $Z' = (Z \times \P V) \cap (S \times \tilde X)$ in $S \times {\Hilb^g} \times \P V$ has the expected dimension,
% Z is flat over Hilb^g, so Z' is flat over X~.
so in $p^* \I_Z$ there are no higher derived pullbacks, only $\I_{Z'/S \times \tilde X}$.
\item For the third term, it is clear that the restriction of $\I_{Z'/C'}$ to $S \times (\zeta,C) \subset S \times \tilde X$ is exactly the sheaf on $S$ parametrized parametrized by $p(\zeta,C) \in \Picbar^{-g}$; thus  we have $\I_{Z'/\C'} = q^* \F \otimes \L$ for some line bundle $\L$ on $\tilde X$, by definition of a moduli space.  Knowing the Picard group of a blow-up, we can write $\L = \O_{\tilde X}(kE) \otimes q^* \L'$ for some $k \in \Z$; then $\L'$ can be absorbed into $\F$, so it is enough to show that $k=1$.

One way to see this is as follows.  Fix a smooth point $\ell \in P \subset \Picbar^{-g}$, and let $\P^1 = q^{-1}(\ell)$.  Then we have $\I_{Z'/\C'}|_{S \times \P^1} = \ell \boxtimes \O_{\P^1}(-k)$.  Observe that $\I_{Z'/\C'}$ has a natural map to $\O_{\C'}$; that $\P^1 = q^{-1}(\ell)$ is naturally identified with $\P \Hom(\ell, \O_C)$, where $C = \operatorname{supp}(\ell) \subset S$; and that on $S \times \P \Hom(\ell, \O_C)$ the natural map must go $\ell \boxtimes \O_{\P^1}(-1) \to \O_C \boxtimes \O_{\P^1}$.  Thus $k=1$ as desired.
% Another way to see it:  Take the first two terms of \eqref{main_seq}, restrict to $S \times \P^1$, tensor with $\O_S(1)$, push down to $\P^1$, and take the determinant.  We have $\chi(\ell(1)) = -1$, so this is $\O_{\P^1}(k)$.  We also calculate that it's $(\det W)^{-1}|_{\P^1}$.  On the other hand we have $\O_{\tilde X}(E) = p^*(\det W) \otimes r^*(\O_{\P V}(-1))$ by adjunction, and restricting this to $\P^1$ we find that $(\det W)|_{\P^1} = \O_{\P^1}(-1)$.  Thus $k=1$.
\end{enumerate}

\noindent Now tensor \eqref{main_seq} with $\O_{\tilde X}(E)$ and apply $p_*$.  The second term becomes $\I_Z$ by the projection formula and the fact that $p_* \O_{\tilde X}(E) = \O_{\Hilb}$.  For the first term, by Grothendieck duality we have
\[ p_*(r^* \O_{\P V}(-1) \otimes \O_{\tilde X}(E))
= (p_* r^* \O_{\P V}(1))^*. \]
From the Koszul resolution of $\tilde X$ in ${\Hilb^g} \times \P V$
\[ 0 \to \det W^* \boxtimes \O(-g) \to \dotsb \to W^* \boxtimes \O(-1) \to \O_{{\Hilb} \times \P V} \to \O_{\tilde X} \to 0 \]
we see that
\[ p_* r^* \O_{\P V}(1) = \cone(W^* \to \O_{\Hilb} \otimes V^*), \]
so from \eqref{underlined_term_seq} we see that its dual is exactly \eqref{underlined_term}, so we have an exact triangle
\begin{equation} \label{second_cone}
\O_S(-1) \boxtimes \pi_{\Hilb,*}(\I_Z \otimes \pi_S^* \O_S(1)) \to \I_Z \to p_*(q^* \F (2E)).
\end{equation}
We would like to conclude $p_*(q^* \F(2E))$ is isomorphic to the cone \eqref{cone_we're_taking}; certainly the two are cones on two non-zero maps between the same pair of objects, so it is enough to argue that $\Hom$ between those objects is 1-dimensional.

The first map of \eqref{second_cone} is obtained by taking the natural map
\[ \O_S(-1) \boxtimes \O_S(1) \to \O_\Delta \]
on $S \times S$ and convolving in the sense of Fourier--Mukai kernels with $\I_Z \in D^b(S \times \Hilb^g)$.  Let $F'\colon D^b(S) \to D^b(\Hilb^g)$ be the $\P^{g-1}$-functor induced by $\I_Z$; then we have
\begin{align*}
&\Hom_{S \times \Hilb}(F' \circ(\O_S(-1) \boxtimes \O_S(1)),\ F') \\
&=\Hom_{S \times S}(\O_S(1) \boxtimes \O_S(-1),\ R' \circ F') \\
&=\Hom_{S \times S}(\O_S(1) \boxtimes \O_S(-1),\ \O_\Delta \oplus \O_\Delta[-2] \oplus \dotsb \oplus \O_\Delta[-2g+2]) \\
&= \CC
\end{align*}
as desired.
\end{proof}

\begin{prop} \label{intertwinement_prop}
If $g \le 5$, so the birational map $p \circ q^{-1}\colon \Picbar^{-g} \dashrightarrow \Hilb^g$ is a Mukai flop, then the diagram
\[ \xymatrix{
D^b(S) \ar[r]^-F \ar[d]_{T_{\O_S(1)}} & D^b(\Picbar^{-g}) \ar[d]^{\KN_2} \\
D^b(S) \ar[r]^-{F'} & D^b(\Hilb^g)
} \]
commutes, where $F$ is induced by $\F$ with the same normalization as in the previous proposition, $F'$ is induced by $\I_Z$, and $\KN_2$ is the equivalence of Definition \ref{KN_k}.\footnote{The reader may be surprised to see $T_{\O_S(1)}$ in this diagram when $T_{\O_S(-1)}$ appeared in \eqref{line_with_T_O(-1)}; the reason for this is as follows.  While we have let $\F$ and $\I_Z$ induce $\P$-functors $F\colon D^b(S) \to D^b(\Picbar^{-g})$ and $F'\colon D^b(S) \to D^b(\Hilb^g)$, from a moduli standpoint it is more natural to let them induce functors $^t\!F\colon D^b(\Picbar^{-g}) \to D^b(S)$ and $^t\!F'\colon D^b(\Hilb^g) \to D^b(S)$, which take skyscraper sheaves of points to the the sheaves on $S$ that they parametrize.  In this direction the cone \eqref{cone_we're_taking} that we studied in Proposition \ref{compat_prop} induces $T_{\O_S(-1)} \circ {}^t\!F'$ as expected.  Transposing to the direction we want, it induces $F' \circ {}^tT_{\O_S(-1)} = F' \circ T_{\O_S(1)}$.}
Thus the $\P$-twists associated to $F$ and $F'$ are conjugate:
\[ P_F' = \KN_2 \circ P_F \circ \KN_2^{-1}. \]
\end{prop}
\begin{proof}
Recall that $\KN_2$ is induced by the line bundle on
\[ \hat X = \tilde X \cup_E (P \times_B P^*) \]
which is $\O(2E)$ on the first component and $\O(-2,-2)$ on the second.  Since the ideal sheaf of $\tilde X$ in $\hat X$ is $\O_{P \times_B P^*}(-E) = \O_{P \times_B P^*}(-1,-1)$, we get an exact sequence of kernels
\begin{equation*} % \label{KN_triangle} if we put in the second half of the remark below.
0 \to \O_{P \times_B P^*}(-3,-3) \to \KN_2 \to \O_{\tilde X}(2E) \to 0.
\end{equation*}
The previous proposition implies that the third term composed with $F$ is $F' \circ T_{\O_S(1)}$, so it remains to show that the first term composed with $F$ is zero.  For this it is enough to restrict $\F$ from $S \times \Picbar^{-g}$ to $S \times P$, tensor with $\O_S \boxtimes \O_P(-3)$, push down to $S \times B$, and show that the result is zero.

By Grothendieck duality we have $\omega_{P/B} = q_* \omega_{E/B}[1]$, so $\O_P(-3)$ differs from $q_* \O_E(-2,-2)[1]$ at most by a line bundle pulled back from $B$, which does not affect this calculation.  Thus we have
\begin{align*}
\F \otimes \O_P(-3)
&\approx \F \otimes q_* \O_E(-2,-2)[1] \\
&= q_*(q^* \F \otimes \O_E(-2,-2))[1] \\
&= q_*(q^* \F(E) \otimes \O_E(-1,-1))[1].
\end{align*}
Tensoring \eqref{main_seq} with $\O_E(-1,-1)$ we get
\opt{ams}{
\begin{multline*}
0 \to \O_S(-1) \boxtimes \O_E(-2,-1) \to p^* \I_Z \otimes \O_E(-1,-1) \\
\to q^* \F(E) \otimes \O_E(-1,-1)\to 0.
\end{multline*}}
\opt{lms}{
\[
0 \to \O_S(-1) \boxtimes \O_E(-2,-1) \to p^* \I_Z \otimes \O_E(-1,-1) \to q^* \F(E) \otimes \O_E(-1,-1)\to 0.
\]}
When we push down to $S \times B$, the first two terms vanish via $p_*$, so the third vanishes as well.
\end{proof}

\begin{rmk}
The end of the last proof implies that the composition
\begin{equation} \label{thing_above}
D^b(S) \xrightarrow{F} D^b(\Picbar^{-g}) \xrightarrow{\varpi_*(\omega_{P/B} \otimes j^*(-))[2]} D^b(B)
\end{equation}
is zero, where $j$ and $\varpi$ are as shown:
\[ \xymatrix{
P \ar[d]_\varpi \ar@{^(->}[r]^j & \Picbar^{-g} \\
{\phantom.}B.
} \]
The second map of \eqref{thing_above} is the right adjoint to the $\P^2$-functor
\[ j_* \varpi^*\colon D^b(B) \to D^b(\Picbar^{-g}), \]
so the image of $j_* \varpi^*$ is orthogonal to the image of $F\colon D^b(S) \to D^b(\Picbar^{-g})$, and thus the $\P^2$-twist around the former commutes with the $\P^{g-1}$-twist around the latter.

%In addition to \eqref{KN_triangle} we have an exact sequence of kernels
%\[ 0 \to \O_{\tilde X}(2E) \to \KN_3 \to \O_{P \times_B P^*}(-3,-3) \to 0, \]
%so the proof of the previous proposition also gives $\KN_3 \circ F = F' \circ T_{\O_S(1)}$.  This agrees with the previous paragraph: the $\P$-twist $P_0$ associated to $j_* \varpi^*$ satisfies $P_0 F = F$, and by Theorem \ref{flop_flop_thm} we have $P_0 = \KN_2^{-1} \circ \KN_3$.
\end{rmk}

% !TeX root = main.tex

\section{Movable cone of the Hilbert scheme} \label{movable}

In the previous section we studied the moduli spaces
\begin{align*}
\Hilb^g &:= \Hilb^g(S) = \M(1,0,1-g) \\
\Picbar^{-g} &:= \Picbar^{-g}(\C/\P V) = \M(0,1,1-2g),
\end{align*}
and the flop $\Hilb^g \dashrightarrow \Picbar^{-g}$ implemented by the spherical twist around $\O_S(-1)$.  In this final section we analyze the movable cone of $\Hilb^g$ using techniques developed by Bayer and Macr\`i in \cite{bm12} and \cite{bm13} and find that these are the only two smooth $K$-trivial birational models.

We recall the description of the Neron--Severi lattice $\NS(\Hilb^g)$ in terms of the Mukai lattice of $S$.  Let
\[ \v_1 := (1,0,1-g) \in H^*_\text{alg}(S,\Z) = H^0(S) \oplus \NS(S) \oplus H^4(S) \]
be the Mukai vector of the ideal sheaf of $g$ points; then the universal ideal sheaf $\I = \I_Z$ on $S \times \Hilb^g$ induces a so-called \emph{Mukai morphism}
\[ \theta_\I\colon \v_1^\perp \to \NS(\Hilb^g), \]
which is an isometry with respect to the Mukai pairing on $\v_1^\perp$ and the Beauville--Bogomolov pairing on $\NS(\Hilb^g)$.  Following \cite{bm12} and \cite{bm13} we adopt the following basis for $\NS(\Hilb^g)$:
\begin{align*}
\tilde H &:= \theta_\I(0,-1,0) & B &:= \theta_\I(-1,0,1-g).
\end{align*}
The divisor $\tilde H$ induces the Hilbert--Chow morphism $\Hilb^g(S) \to \operatorname{Sym}^g(S)$, and $B$ is half of the exceptional divisor of that morphism.

The flop $\Hilb^g \dashrightarrow \Pic^{-g}$ is regular away from a codimension-2 locus, so it lets us identify $\NS(\Picbar^{-g})$ with $\NS(\Hilb^g)$ and $\Mov(\Picbar^{-g})$ with $\Mov(\Hilb^g)$. 

\begin{prop}
The movable cone of $\Hilb^g$ is generated by $\tilde H$ and $\tilde H - B$, and is divided into two chambers as follows: 
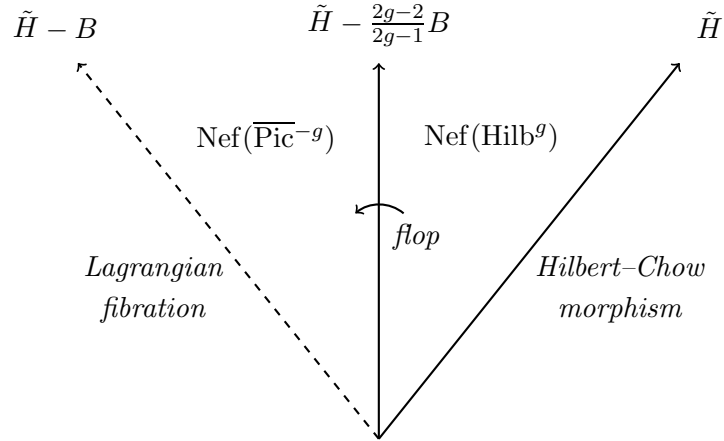
\vspace{1ex} % for some reason having \\ instead of \vspace here gave an underfull hbox
\begin{center}
\begin{tikzpicture}[scale=1]
\draw [->][dashed,thick] (0,0) -- (-4,5);
\draw [->][thick] (0,0) -- (0,5);
\draw [->][thick] (0,0) -- (4,5);
\node at (1.5,4) {$\Nef(\Hilb^g)$};
\node at (-1.5,4) {$\Nef(\Picbar^{-g})$};
\node at (4.4,5.5) {$\tilde H$};
\node at (0,5.5) {$\tilde H - \frac{2g-2}{2g-1}B$};
\node at (-4.3,5.5) {$\tilde H - B$};
\node at (3.2,2) {\begin{tabular}{c} Hilbert--Chow \\ morphism\end{tabular}}; 
\draw [->](0.33,3)[thick] arc (50:130:0.5);
\node[anchor=north] at (0.5,3) {flop};
\node at (-3,2) {\begin{tabular}{c} Lagrangian \\ fibration\end{tabular}};
\end{tikzpicture}
\end{center}
Thus $\Hilb^g$ has no other smooth $K$-trivial birational models apart from $\Picbar^{-g}$.
\end{prop}
\begin{proof}
The generators of $\Mov(\Hilb^g)$ can be read off from \cite[Prop.~13.1(a)]{bm13}.  The generators of $\Nef(\Hilb^g)$ can be read off from \cite[Prop.~10.3]{bm12}.\footnote{Note that when $g=2$ this seems to disagree with \cite[Lem.~13.3(b)]{bm13}, but the latter has a typo, as is clear from the example that follows it.  It should read ``Otherwise, let $(x_1,y_1)$ be the minimal positive solution of (34).  Then $\Nef(M) = \langle \tilde H, \tilde H - 2d\frac{y_1}{x_1} B \rangle$.''  Or else equation (34) should read ``$X^2 - dY^2 = 5$.''}  The pullback of a hyperplane via the Lagrangian fibration $\Picbar^{-g} \to \P V$ must be one wall of $\Nef(\Picbar^{-g}$); moreover it is isotropic for the Beauville--Bogomolov form, and the only isotropic ray in $\Mov(\Hilb^g)$ is $\tilde H - B$.  Thus it remains to show that $\Nef(\Picbar^{-g})$ extends all the way to $\tilde H - \frac{2g-2}{2g-1} B$.

On the one hand, by \cite[Thm.~12.1]{bm13} this is a purely lattice-theoretic question: walls in $\Mov(\Hilb^g)$ arise from vectors $\a \in H^*_\text{alg}(S,\Z)$ satisfying certain numerical conditions.  For each fixed $g$ it is straightforward to show that there are no more such $\a$, but we do not see how to do it for all $g$ at once.

Instead we use \cite[Example~9.7]{bm12} to produce a one-parameter family of ample divisors on $\Picbar^{-g}$ that fill out the cone from $\tilde H - B$ to $\tilde H - \frac{2g-2}{2g-1} B$.  Let $\v_2 = (0,1,1-2g)$, let $\F$ be a universal sheaf on $S \times \Picbar^{-g}$, and let
\[ \theta_\F\colon \v_2^\perp \to \NS(\Picbar^{-g}) \]
be the associated Mukai morphism.  Apply \cite[Example~9.7]{bm12} with $d=g-1$, $v = \v_2$, $A = 1$, and $B = -1$; then we find ourselves in ``Case 2,'' and we get $U = \O_S(-1)$ and $t_0 = \frac{1}{\sqrt{g-1}}$.  The conclusion is that $\theta_\F(w_{\sigma_{t,-1}})$ is ample for all $t > t_0$, where $w_{\sigma_{t,-1}}$ can be calculated using [\emph{ibid.}, Lem.~9.2]:
\[ w_{\sigma_{t,-1}} = (2g-2)t \cdot \left(
1,\,
-\tfrac{2g-1}{2g-2},\,
g - (g-1)t^2
\right). \]
In order to express $\theta_\F(w_{\sigma_{t,-1}})$ in terms of $\tilde H$ and $B$, we must compare $\theta_\F$ to $\theta_\I$.  On the open set where the flop $\Hilb^g \dashrightarrow \Picbar^{-g}$ is defined, the spherical twist around $\O_S(-1)$ turns $\I$ into $\F$; thus the reflection through the hyperplane orthogonal to $\s := v(\O_S(-1)) = (1,-1,g)$ turns $\theta_\I$ into $\theta_\F$ and vice versa.  Recalling that this reflection is $\v \mapsto \v + \langle \v, \s \rangle\, \s$, we calculate that our ample divisor $\theta_\F(w_{\sigma_{t,-1}})$ is, up to rescaling by a positive number,
\[ \theta_\I\left(1,\, -1 - \tfrac{1}{2(g-1)^2 t^2},\, g-1 \right)
= (1 + \tfrac{1}{2(g-1)^2 t^2}) \tilde H - B. \]
Letting $t \to t_0 = \frac{1}{\sqrt{g-1}}$ we get a nef divisor $\tfrac{2g-1}{2g-2} \tilde H - B$, which is a positive multiple of $\tilde H - \tfrac{2g-2}{2g-1} B$ as claimed.  Letting $t \to \infty$ we rediscover the Lagrangian fibration wall $\tilde H - B$.
\end{proof}

\begin{rmk}
The nef cones of $\Hilb^g$ and $\Picbar^{-g}$ can also be obtained from a result of Yoshioka \cite[Prop.~4.39]{yoshioka}.  Note that his Mukai morphism differs by a sign from the one used here.
\end{rmk}

\begin{rmk}
In general, the number of smooth $K$-trivial birational models of $\Hilb^n(S)$ as $n$ varies is interesting and hard to control.  For example, it follows from \cite[Lem.~13.3]{bm13} that $\Hilb^2(S)$ has more than one model for $g = 2, 6, 12$ but only one model for other $g < 20$.  On the other hand, if $g=2$ then $\Hilb^7(S)$ has \emph{six} different models, related by various complicated flops \cite[Example~13.5]{bm13}.
\end{rmk}

\newcommand \httpurl [1] {\href{http://#1}{\nolinkurl{#1}}}
\bibliographystyle{plain}
\bibliography{main}

\opt{lms}{
\affiliationone{
Nicolas Addington \\
Department of Mathematics \\
University of Oregon \\
Eugene, OR 97403-1222 \\
United States
\email{adding@uoregon.edu}}
\affiliationtwo{
Will Donovan \\
Kavli Institute for the Physics and Mathematics of the Universe (WPI) \\
The University of Tokyo Institutes for Advanced Study \\
5-1-5 Kashiwanoha \\
Kashiwa, Chiba, 277-8583 \\
Japan
\email{will.donovan@ipmu.jp}}
\affiliationthree{
Ciaran Meachan \\
School of Mathematics \\
The University of Edinburgh \\
James Clerk Maxwell Building \\
Peter Guthrie Tait Road \\
Edinburgh EH9 3FD \\
United Kingdom
\email{ciaran.meachan@ed.ac.uk}}
}

\end{document}